\documentclass[a4, paper]{article}

\usepackage{amsmath,amsfonts}

\usepackage[T1]{fontenc}

\title { $ \sup + \inf $ for Riemannian surfaces and $ \sup \times \inf $ for bounded domains of $ {\mathbb R}^n $, $ n\geq 3 $  }

\author {Samy Skander Bahoura \footnote {email: samybahoura@yahoo.fr, bahoura@ccr.jussieu.fr}}

\date {\small { Universit\'e Paris VI, 4 place Jussieu 75005 Paris, France}}

\begin{document}

\maketitle

\rule{120mm}{0.5mm}

\smallskip

\underbar { Abstract}

On a Riemannian surface, we give a condition to obtain a minoration of $ \sup +\inf $. On an open bounded set of $ {\mathbb R}^n $ $ (n\geq 3) $ with smooth boundary, we have a minoration of $\sup \times \inf $ for prescribed scalar curvature equation with Dirichlet condition.

\smallskip


\rule {120mm}{0.5mm}

\bigskip

In this paper, we study some inequalities of type  $ \sup +\inf $ (in dimension 2) and $ \sup \times \inf $ (in dimension $ n\geq 3 $). We denote $ \Delta =-\nabla_i(\nabla^i) $ the geometric laplacian.

\smallskip

The paper is linking to the Note presented in Comptes Rendus de l'Acad\'emie des Sciences de Paris (see [B1])

\smallskip

In dimension 2, we work on Riemannian surface $ (M,g) $ and we consider the following equation:

$$ \Delta u+f=Ve^u \qquad (E_1) $$

where $ f $ and $ V $ are two functions.

We are going to prove a minoration of $ \sup u+\inf u $ under some conditions on  $ f $ and $ V $.

\bigskip

Where $ f=R $, with $ R $ the scalar curvature of $ M $, we have the scalar curvature equation studied by T. Aubin, H. Brezis, YY. Li, L. Nirenberg, R. Schoen.

\smallskip

In the case $ f= R = 2 \pi $ and $ M={\mathbb S}^2 $, we have a lower bound for $ \sup + \inf $ assuming $ V $ non negative, bounded above by a positive constant $ b $ and without condition on $ \nabla V $ (see  Bahoura [B]).

\smallskip

The problem was studied when we suppose $ V=V_i $ uniformly lipschitzian and between two positive constants. (See Bahoura [B] and Li [L]). In fact, there exists  $ c=c(a,b,A,M) $ such that for all sequences $ u_i $ and $ V_i $ satisfying:

$$  \Delta u_i+R=V_ie^{u_i},\,\, 0< a\leq V_i(x)\leq b \,\, {\rm and}\,\, ||\nabla V_i||_{\infty} \leq A ,$$

we have,

$$ \sup_M u_i+\inf_M u_i \geq c \,\, \forall \,\, i.$$

We have some results about $ L^{\infty} $ boundness and asymptotic behavior for the solutions of euqations of this type on open set of $ {\mathbb R}^2 $, see [BM], [S], [SN 1] and [SN 2].

\bigskip

Here, we try to study the same problem with minimal conditions on $ f $ and $ V $, we suppose $ 0 \leq V \leq b  $ and without assumption on $ \nabla V $.

\bigskip

{\it Theorem 1. Assume $ (M,g) $ a Riemannian surface and $ f,V $ two functions satisfying:

$$ f(x)\geq 0,\,\,{\rm and}\,\, 0\leq V(x) \leq b <+\infty,\,\,\, \forall \,\, x\in M . $$

suppose $ u $ solution of:

$$ \Delta u+f=Ve^u . $$

then:

if $ 0< \int_M f \leq 8\pi $, there exists a constant $ c=c(b,f,M) $ such that:

$$ \sup_M u+\inf_M u \geq c , $$

if $ 8 \pi <\int_M f < 16 \pi $, there exists $ C=C(f,M) \in ]0,1[ $ and $ c=c(b,f,M) $ such that:

$$ \sup_M u+C\inf_M u \geq c .$$}

{\bf Remark:} In fact, we can suppose $ f\equiv k $ a constant. (See [B1]).

\bigskip

Now, we work on a smooth bounded domain $ \Omega \subset {\mathbb R}^n $ $ (n\geq 3) $.

\smallskip

Let us consider the following equations:

$$ \Delta u_{\epsilon}={u_{\epsilon}}^{N-1-\epsilon},\,\,\, u_{\epsilon}>0\,\,\,{\rm in} \,\,\,\Omega \,\,\, {\rm and} \,\,\, u_{\epsilon}=0 \,\,\, {\rm on } \,\, \partial \Omega \qquad  (E_2). $$

with $ \epsilon \geq 0 $, $ N=\dfrac{2n}{n-2} $.

\smallskip

The existence result for those equations depends on the geometry of the domain. For example, if we suppose, $ \Omega $ starshaped and $ \epsilon =0 $, the Pohozaev identity assure a nonexistence result. If $ \epsilon =0 $, under assumption on $ \Omega $, we can have an existence result. When $ \epsilon >0 $ there exists a solutions for the previous equation.

\bigskip

For $ \epsilon >0 $,  [AP], [BP] and [H] , studied some properties of the previous equation.

\bigskip 

On unit ball of $ {\mathbb R}^n $, Atkinson-Peletier(see [AP]) have proved:

$$ \lim_{\epsilon \to 0} [\sup_{B_1(0)} u_{\epsilon} \inf_{B_k(0)} u_{\epsilon} ]=  \left (\dfrac{1}{|k|}-1 \right ) ,$$

with $ |x|=k<1 $.

\bigskip

In [H], Z-C Han, has proved the same estimation on a smooth open set $ \Omega \subset {\mathbb R}^n $ with the following condition :

$$ \lim_{\epsilon \to 0} \dfrac{\int_{\Omega} |\nabla u_{\epsilon}|^2}{[||u_{\epsilon}||_{L^{N-1-\epsilon}}]^2}=S_n \qquad \qquad (1), $$

with $ S_n=\pi n(n-2)\left [\dfrac{\Gamma(n/2)}{\Gamma (n)} \right ] $ the best constant in the Sobolev imbedding.

\bigskip

In fact, the result of Z-C Han (see [H] ), is (with (1)),

$$ \lim_{\epsilon \to 0}||u_{\epsilon}||_{L^{\infty}} u_{\epsilon}(x)= \sigma_n (n-2) G(x,x_0),\,\,\, {\rm with },\,\, x\in \Omega-\{x_0\} .$$

where $ x_0 \in \Omega $ and $ G $ is the Green function with Dirichlet condition.

\smallskip

In our work, we search to know if it is possible to have a lower bound of $ \sup \times \inf $, without the assumption $ (1) $.

\bigskip

{\it Theorem 2. For all compact $ K $ of $ \Omega $, there exists a positive constant  $ c=c(K,\Omega, n)>0 $, such that for all solution $ u_{\epsilon} $ of $ (E_2) $ with $ \epsilon \in ]0,\dfrac{2}{n-2} ] $, we have: 

$$ \sup_{\Omega} u_{\epsilon} \times \inf_K u_{\epsilon} \geq c. $$}

Next, we are intersting by the following equation:

$$ \Delta u=u^{N-1}+\epsilon u,\,\,\, u>0, \,\,\, {\rm in }\,\,\, \Omega,\,\,\, {\rm and}\,\,\, u=0 \,\,\,{\rm on} \,\,\, \partial \Omega. $$

We know that in dimension 3, there is no radial solution for the previous equation if $ \epsilon \leq \lambda_* $ with $ \lambda_*>0 $, see [B N]. Next, we consider $ n\geq 4 $.

\bigskip

We set $ G $ the Green function of the laplacian with Dirichlet condition. For $ 0< \alpha <1 $, we denote:

$$ \beta=\dfrac{\alpha}{\sup_{\Omega} \int_{\Omega} G(x,y)dy} .$$

Assume $ n\geq 4 $, we have:

\bigskip

{\it Theorem 3. For all compact $ K $ of $ \Omega $ and all $ 0< \alpha <1 $  there is a positive constant $ c=c(\alpha, K,\Omega, n)  $ , such that for all sequences $ (\epsilon_i)_{i\in {\mathbb N}} $ with $ 0 < \epsilon_i \leq \beta $ and $ (u_{\epsilon_i})_{i\in {\mathbb N}}) $ satisfying:

$$  \Delta u_{\epsilon_i}={u_{\epsilon_i}}^{N-1}+\epsilon_i u_{\epsilon_i}, \,\,\, u_{\epsilon_i}>0\,\,\, {\rm and }\,\,\, u_{\epsilon_i}=0 \,\,\,{\rm on} \,\,\, \partial \Omega,\,\,\, \forall \,\,\, i ,$$}

we have:

$$ \forall \,\,\, i,\,\,\, \sup_{\Omega} u_{\epsilon_i} \times \inf_K u_{\epsilon_i} \geq c .$$

\newpage

\underbar{{\bf Proof of the Theorem 1:}}

\smallskip

\underbar {\bf First part ( $ 0<\int_M f \leq 8\pi $ ) :}

\smallskip

We have:

$$ \Delta u+f=Ve^u, $$

We multiply by $ u $ the previous equation and we integrate by part, we obtain:

$$ \int_M |\nabla u|^2+\int_M fu=\int_M Ve^u u, $$

But $ V\geq 0 $ and $ f\geq 0 $, then :

$$ \int_M |\nabla u|^2+\inf_M u \int_M f \leq \sup_M u \int_M V e^u. $$

On Riemannian surface, we have the following Sobolev inequality, (see [DJLW], [F]):

$$ \exists\,\, C=C(M,g) >0, \,\,\forall v\in H_1^2(M),\,\, \log \left ( \int_M e^v \right )\leq \dfrac{1}{16\pi}\int_M |\nabla v|^2+\dfrac{1}{Vol(M)} \int_M v +\log C .$$

Let us consider $ G $ the Green function of the laplacian such that:

 $$ G(x,y)\geq 0 \,\,\,\,  {\rm and },\,\, \int_M G(x,y)dV_g(y) \equiv k = {\rm constant} .$$
 
 Then,
 
$$  u(x)=\dfrac{1}{Vol(M)} \int_M u+\int_M G(x,y)[V(y)e^{u(y)}-f(y)] dV_g(y),$$

and,

$$ \inf_M u=u(x_0)\geq \dfrac{1}{Vol(M)}\int_M u-C_1, $$

with,

$$\int_M [G(x_0,y)f(y)]\leq \sup_M f \int_M G(x_0,y)dV_g(y)=k \sup_M f= C_1 .$$

But, $ \int_M V e^u =\int_M f >0 $, we obtain,

$$ \left( \int_M f \right)(\sup_M u+\inf_M u)\geq -2C_1\int_M f+  \dfrac{2}{Vol(M)}\left (\int_M u \right)\left (\int_M f\right)+\int_M |\nabla u|^2, $$

thus,

$$ \sup_M u+\inf_M u\geq 2\left [\dfrac{1}{Vol(M)} \int_M u+\dfrac{1}{2\int_M f} \int_M |\nabla u|^2 \right]-2C_1. $$

If we suppose, $ 0< \int_M f\leq 8\pi $, we obtain $ \dfrac{1}{2\int_M u} \geq \dfrac{1}{16 \pi} $  and then:

$$ \sup_M u+\inf_M u \geq 2\left [ \dfrac{1}{Vol(M)} \int_M u +\dfrac{1}{16\pi} \int_M |\nabla u|^2\right ]-2C_1, $$

We use the previous Sobolev inequality, we have:

$$ \sup_M u+\inf_M u \geq -2C_1-2\log C+2\log \left (\int_M e^u \right ),$$

but,

$$ \int_M f=\int_M Ve^u \leq b\int_M e^u, $$

then,

$$ \int_M e^u \geq \dfrac{1}{b}\int_M f, $$

and finaly,

$$ \sup_M u+\inf_M u\geq -2C_1-2\log C+2\log \left (\dfrac{1}{b} \int_M f \right ) .$$

\underbar {\bf Second part ( $ 8\pi< \int_M f < 16\pi $ ):}

\smallskip

Like en the first part, we have:

$$ a) \,\,\, \int_M |\nabla u|^2+\inf_M u\int_M f \leq \sup u \int_M f ,$$

$$ b)\,\,\, \log \left ( \int_M e^u \right )\leq \dfrac{1}{16 \pi} \int_M |\nabla u|^2+\dfrac{1}{Vol(M)} \int_M u+\log C , $$

$$ c) \,\,\, \inf_M u\geq \dfrac{1}{Vol(M)} \int_M u -C_1 .$$

We set $ \lambda >0 $. We use $ a), b), c) $ and we obtain:

$$ \left ( \int_M f \right )( \sup_M u+\lambda \inf_M u ) \geq -(\lambda +1) C_1 \int_M f+ \dfrac{(1+\lambda )}{Vol(M)} \left ( \int_M u \right )\left (\int_M f \right )+\int_M |\nabla u|^2, $$

thus,

$$ \sup_M u+\lambda \inf_M u \geq -(\lambda+1) C_1+(1+\lambda ) \left [\dfrac{1}{Vol(M)} \int_M u+ \dfrac{1}{(1+\lambda) \int_M f} \int_M |\nabla u|^2 \right ].$$

We choose  $\lambda>0 $, such that, $ \dfrac{1}{(1+\lambda) \int_M f} \geq \dfrac{1}{16 \pi}, $

thus, $ (1 +\lambda)\int_M f \leq 16 \pi $, $ 0 < \lambda \leq \dfrac{16\pi-\int_M f}{\int_M f} < 1 $.

Finaly, the choice of $ \lambda $, give:

$$ \sup_M u+\lambda \inf_M u \geq -(\lambda+1)C_1-(1+\lambda)\log C+(1+\lambda ) \log \left (\dfrac{1}{b}\int_M f\right ). $$

If we take $ \lambda =\dfrac{16 \pi-\int_M f}{\int_M f}\in ]0,1[ $, we obtain:

$$ \sup_M u+\left ( \dfrac{16\pi-\int_M f}{\int_M f} \right )\inf_M u \geq -C_1 \dfrac{16\pi}{\int_M f}-\dfrac{16\pi}{\int_M f}\log C+\dfrac{16\pi}{\int_M f} \log \left (\dfrac{1}{b} \int_M f \right).$$

\bigskip

\bigskip

\underbar{{ \bf Proof of theorems 2 and 3:}}

\bigskip

Here, we give two methods to prove the theorems 2 and 3, but we do the proof only for the theorem 2. In the first proof we use the Moser iterate scheme, the second proof is direct.

\bigskip

\underbar {\bf Method 1:} by the Moser iterate scheme.

\bigskip

We argue by contradiction and we suppose:

$ \exists \,\, K \subset \, \subset \,\Omega,\,\, \forall \, c>0,\,\, \exists\,\epsilon_c \in ]0, \dfrac{2}{n-2}] $ such that:

$$ \Delta u_{{\epsilon}_c}={u_{{\epsilon}_c}}^{N-1-{\epsilon}_c},\,\, u_{\epsilon_c} >0 \,\, {\rm in } \,\, \Omega \,\, {\rm and} \,\,  u_{\epsilon_c}=0 \,\, {\rm on }\,\, \partial\Omega,  $$

with,

$$ \sup_{\Omega} u_{{\epsilon_c}} \times \inf_K u_{{\epsilon}_c} \leq c $$  

We take $ c=\dfrac{1}{i} $, there exists a sequence $ (\epsilon_i)_{i\geq 0} $, such that $ \forall \,\, i\in {\mathbb N}, \epsilon_i \in ]0,\dfrac{n}{n-2}]$ and

$$ \Delta u_{\epsilon_i} ={u_{\epsilon_i}}^{N-1-\epsilon_i}, \,\, u_{\epsilon_i} >0\,\,{\rm in } \,\, \Omega\,\, {\rm and} \,\, u_{\epsilon_i}=0, \,\, {\rm on} \,\,\partial  \Omega \,\, (*) $$

with, 

$$ \sup_{\Omega} u_{\epsilon_i} \times \inf_K u_{\epsilon_i} \leq \dfrac{1}{i} \to 0 \,\,(**) .$$

Clearly the function $ u_{\epsilon_i} $ which satisfy $(*)$, there exists $ x_{\epsilon_i} \in \Omega $ such that: 

$$ \sup_{\Omega} u_{\epsilon_i} =\max_{\Omega} u_{\epsilon_i}=u_{\epsilon_i}(x_{\epsilon_i}).$$

\underbar{{\bf Lemma:}}

There exists $ \delta=\delta(\Omega,n)>0 $ such that for all $ \epsilon >0 $ and $ u_{\epsilon }>0 $, solution of our problem with $ x_{\epsilon}\in\Omega $, $ \sup_{\Omega} u_{\epsilon}=u_{\epsilon}(x_{\epsilon})$ we have:

$$ d(x_{\epsilon}, \partial \Omega) \geq \delta. $$

\underbar{{\bf Proof of the lemma:}}

\smallskip

We argue by contradiction. We suppose: $ \forall \,\, \delta >0 $, $ \exists $ $ x_{\epsilon_{i_{\delta}}} $ such that: $ d(x_{\epsilon_{i_{\delta}}},\partial \Omega ) \leq \delta $.

\bigskip

We take $ \delta =\dfrac{1}{j}, j\to +\infty $, we have a subsequence $ \epsilon_{i_j} $, noted $ \epsilon_i $, such that, $ d(x_{\epsilon_i},\partial \Omega )\to 0 $ .

\bigskip

Let us consider $ G $ the Green function of the laplacian with Dirichlet condition and $ w $ satisfying:

$$ \Delta w=1\,\, {\rm in} \,\, \Omega \,\, {\rm and} \,\, w=0 \,\, {\rm on } \,\, \partial \Omega .$$

Using the variational method, we can prove the existence of $ w $ and $ w\in {\cal C}^{\infty}({\bar {\Omega}}) $.

\bigskip

The Green representation formula and the fact $ x_{\epsilon_i} \to y_0 \in \partial \Omega $ give:

$$  0=w(y_0) \leftarrow w(x_{\epsilon_i})=\int_{\Omega} G(x_{\epsilon_i},y)dy ,$$

we can write,

$$ \int_{\Omega} G(x_{\epsilon_i},y)dy \to 0. $$

The function $ u_{\epsilon_i} $ satisfy $ (*) $ and thus:

$$ u_{\epsilon_i} (x_{\epsilon_i})\leq (\max_{\Omega} u_{\epsilon_i})^{N-1-\epsilon_i} \int_{\Omega} G(x_{\epsilon_i},y) dy ,$$

consequently,

$$ 1\leq [u_{\epsilon_i}(x_{\epsilon_i})]^{N-2-\epsilon_i} \int_{\Omega} G(x_{\epsilon_i},y)dy. $$

Then,

$$ u_{\epsilon_i}(x_{\epsilon_i})\to +\infty \,\, {\rm and}\,\, x_{\epsilon_i}\to y_0 \in \, \partial \Omega \,\,\, (***). $$

But, if we use the result of Z-C.Han (see [H] page 164) and [DLN] (pages 44-45 and 50-53) and the moving plane method (see [GNN]) we obtain:

\bigskip

if $ \Omega $ is smooth bounded domain, $ f $ a function in $ {\cal C}^1 $ and  $ u $ is a solution of:

$$ \Delta u=f(u),\,\, {\rm in} \,\, \Omega \,\, {\rm and} \,\, u=0 \,\, {\rm on}\,\, \partial \Omega ,$$

there exists two positive constants $ \delta $ and $ \gamma $, which depend only on the geometry of the domain $ \Omega $, such that:

\bigskip

$ \forall \,\, x \in \{ z, d(z, \partial \Omega)\leq \delta \}, \,\, \exists \,\, \Gamma_x \subset \{ z, d(z, \partial \Omega)\geq \dfrac{\delta}{2} \}$ with $ mes(\Gamma_x)\geq \gamma $ et $ u(x)\leq u(\xi) $ for all $ \xi \in \Gamma_x $.

\bigskip

Thus, 

$$ u(x)\leq \dfrac{1}{mes(\Gamma_x)} \int_{\Gamma_x} u \leq  \dfrac{1}{\gamma} \int_{\Omega'} u \qquad (*') , $$
 
with $ \Omega'\subset \subset \Omega $.

\smallskip

If we replace $ x $ by $ x_{\epsilon_i} $, $ u $ by $ u_{\epsilon_i} $ and we take $ \Omega'=\{z\in \Omega, d(z,\partial \Omega)\geq \dfrac{\delta}{2} \} $, we obtain (after using the argument of the first eigenvalue like in [H]):

$$ +\infty \leftarrow u_{\epsilon_i}(x_{\epsilon_i})\leq  \dfrac{1}{\gamma} \int_{\Omega'} u_{\epsilon_i} \leq c_2(\Omega',n)<\infty,$$

it is contradiction. The lemma is proved.

\bigskip

We continue the proof of the Theorem.

\bigskip

Without loss of generality , we can assume $ x_{\epsilon_i} \to y_0 $. We consider $ (x_{\epsilon_i})_{i\geq 0} $ and $ \mu>0 $, such that $ x_{\epsilon_i} \in B(y_0,\mu)\subset \subset \Omega $. ( we take $ \mu=\dfrac{\delta}{2} $ for example). 

\bigskip

We have:

$$ u_{\epsilon_i}(x)=\int_{\Omega} G(x,y){u_{\epsilon_i}}^{N-1-\epsilon_i}(y) dy $$

According to the properties of the Green functions and maximum priciple, on the compact $ K $ of $ \Omega $:

$$ G(x,y)\geq c_3=c(K,\Omega,n)>0,\,\, \forall \,\, x\in K,\,\, y\in B(y_0,\mu).$$

Thus,

$$ \inf_K u_{\epsilon_i}=u_{\epsilon_i}(y_{\epsilon_i}) \geq  c_3 \int_{B(y_0,\mu)} {u_{\epsilon_i}}^{N-1-\epsilon_i},$$

and then,

$$ \int_{B(y_0,\mu)} {u_{\epsilon_i}}^{N-\epsilon_i} \leq  ( \sup_{\Omega} u_{\epsilon_i}) \times  \int_{B(y_0,\mu)} {u_{\epsilon_i}}^{N-1-\epsilon_i} \leq \dfrac{ (\sup_{\Omega} u_{\epsilon_i} \times \inf_K u_{\epsilon_i})}{c_3} \to 0.$$

Finaly,

$$ 0 < \int_{B(y_0,\mu)} {u_{\epsilon_i}}^{N-\epsilon_i} \to 0 \qquad (****).$$

Let $ \eta $ be a smooth function such that :

$$ 0 \leq \eta \leq 1,\,\, \eta \equiv 1,\,\, {\rm on} \,\, B(y_0,\mu/2),\,\, \eta \equiv 0,\,\, {\rm on } \,\, \Omega-B(y_0,\dfrac{2\mu}{3}).$$    

Set $ k>1 $. We multiply the equation of $ u_{\epsilon_i} $ by ${u_{\epsilon_i}}^{2k-1} {\eta}^2 $ and we integrate by part the first member,

$$ (2k-1) \int_{B(y_0,2\mu/3)} |\nabla u_{\epsilon_i}|^2 {u_{\epsilon_i}}^{2k-2}{\eta }^2+2\int_{B(y_0,2\mu/3)} <\nabla u_{\epsilon_i}|\nabla \eta > \eta {u_{\epsilon_i}}^{2k-1} =\int_{B(y_0,2\mu/3)} {u_{\epsilon_i}}^{N-2k-2-\epsilon_i}{\eta }^2 ,$$

We compute $ |\nabla ({u_{\epsilon_i}}^k \eta)|^2 $ and we deduce:

$$ \dfrac{2k-1}{k^2}\int_{B(y_0,2\mu/3)} |\nabla ({u_{\epsilon_i}}^k \eta)|^2+\dfrac{2-2k}{k} \int_{B(y_0,2\mu/3)} <\nabla u_{\epsilon_i}|\nabla \eta > {u_{\epsilon_i}}^{2k-1} \eta-\dfrac{2k-1}{k^2}\int_{B(y_0,2\mu/3)} |\nabla \eta |^2{u_{\epsilon_i}}^{2k} $$

$$ =\int_{B(y_0,2\mu/3)}{\eta }^2 {u_{\epsilon_i}}^{N+2k-2-\epsilon_i} .$$

And,

$$ \int_{B(y_0,2\mu/3)} <\nabla u_{\epsilon_i}|\nabla \eta > {u_{\epsilon_i}}^{2k-1} \eta =\dfrac{1}{4k}\int_{B(y_0,2\mu/3)} <\nabla ({u_{\epsilon_i}}^{2k})|\nabla ({\eta }^2)>=\dfrac{1}{4k} \int_{B(y_0,2\mu/3)} \Delta ({\eta}^2){u_{\epsilon_i}}^{2k} .$$

Then,

$$ \dfrac{2k-1}{k^2} \int_{B(y_0,2\mu/3)}|\nabla ({u_{\epsilon_i}}^{k}\eta )|^2=\dfrac{2-2k}{4k^2}\int_{B(y_0,2\mu/3)} \Delta ({\eta}^2 ){u_{\epsilon_i}}^{2k}+\dfrac{2k-1}{k^2} \int_{B(y_0,2\mu/3)}|\nabla \eta |^2 {u_{\epsilon_i}}^{2k}+ $$

$$ +\int_{B(y_0,2\mu/3)}{u_{\epsilon_i}}^{N+2k-2-\epsilon_i} .$$

But,

$$ \int_{B(y_0,2\mu/3)} {u_{\epsilon_i}}^{N+2k-2-\epsilon_i}=\int_{B(y_0,2\mu/3)} ({u_{\epsilon_i}}^{2k} {\eta }^2)({u_{\epsilon_i}}^{N-2-\epsilon_i}) . $$

Using H\"older inequality with $ p=(N-\epsilon_i)/2 $ and $ p'=(N-\epsilon_i)/(N-\epsilon_i-2) $, we obtain:

$$\dfrac{2k-1}{k^2} [||\nabla (\eta u_{\epsilon_i})||_{L^2(B_0)}]^2 \leq [||u_{\epsilon_i}||_{L^{N-\epsilon_i}(B_0)}]^{N-\epsilon_i-2} \times [ ||\eta {u_{\epsilon_i}}^k||_{L^{N-\epsilon_i}(B_0)}]^2+C [||u_{\epsilon_i}||_{L^{2k}(B_0)}]^{2k} $$

with $ B_0=B(y_0,2\mu/3) $ and $ C=C(k,\eta)=\dfrac{2-2k}{4k^2}||\Delta \eta||_{\infty}+\dfrac{2k-1}{k^2}||\nabla \eta||_{\infty} $.

\bigskip

H\"older and Sobolev inequalities give,

$$ [||\eta{u_{\epsilon_i}}^k||_{L^{N-\epsilon_i}(B_0)}]^2 \leq |B_0|^{2\epsilon_i/[N(N-\epsilon_i)]} K [||\nabla (\eta {u_{\epsilon_i}}^k)||_{L^2(B_0)}]^{2} .$$

We obtain:

$$ \dfrac{2k-1}{K k^2 |B_0|^{2\epsilon_i/[N(N-\epsilon_i)]}} [||\eta {u_{\epsilon_i}}^k||_{L^{N-\epsilon_i}(B_0)}]^{2} \leq [||u_{\epsilon_i}||_{L^{N-\epsilon_i}(B_0)}]^{N-2-\epsilon_i} \times [|| \eta {u_{\epsilon_i}}^k ||_{L^{N-\epsilon_i}(B_0)}]^2+ $$

$$ +C(k,\eta))[||u_{\epsilon_i}||_{L^{2k}(B_0)}]^{2k} , $$

with $ |B_0|=mes[B(0,2\mu /3)] $. 

\bigskip

We choose $ k=\dfrac{N-\epsilon_i}{2} $ and we denote $ \alpha_i=[||\eta {u_{\epsilon_i}}^{(N-\epsilon_i)/2}||_{L^{N-\epsilon_i}(B_0)}]^2 >0 $.

\bigskip

We have:

$$ c_1 \alpha_i\leq \beta_i\alpha_i+c_2\gamma_i, $$

with $ c_1=c_1(N,\mu)>0, c_2=c_2(N,\mu)>0 $, $ \beta_i =[||u_{\epsilon_i}||_{L^{N-\epsilon_i}}]^{N-2-\epsilon_i} $ and  $ \gamma_i=[||u_{\epsilon_i}||_{L^{N-\epsilon_i}}]^{N-\epsilon_i} $.

\bigskip

with $ \epsilon_i\in ]0,\dfrac{2}{n-2}] $. According to $ (****) $, we have, $\beta_i\to 0 $ and $ \gamma_i\to 0 $.

\bigskip

Thus,

\bigskip

$$ (c_1/2)\alpha_i \leq (c_1-\beta_i)\alpha_i \leq \gamma_i \to 0 .$$

Finaly, 

$$ 0 < \int_{B(y_0,\mu/2)} {u_{\epsilon_i}}^{(N-\epsilon_i)^2/2} \leq \int_{B(y_0,2\mu/3)} {\eta u_{\epsilon_i}}^{(N-\epsilon_i)^2/2} \to 0 .$$

We iterate this process with $ k=\dfrac{(N-\epsilon_i)^2}{4} $ after with $ k=\dfrac{(N-\epsilon_i)^r}{2^r} $, $ r\in {\mathbb N}^* $, we obtain, for all $ q \geq 1 $, there exists $ l>0 $, such that:

$$ \int_{B(y_0,l)} (u_{\epsilon_i})^q \to 0.$$

Using the Green representation formula, we obtain:

$$ \forall \,\, x\in B(x,l'),\,\, u_{\epsilon_i}(x)=\int_{B(y_0,l)} G(x,y){u_{\epsilon_i}}^{N-1-\epsilon_i}(y)dy+\int_{\partial B(y_0,l)} \partial_{\nu} G(x,\sigma_l) {u_{\epsilon_i}}(\sigma_l) d\sigma_l \,\,\, (*****).$$

where $ 0 < l'\leq l $.

\smallskip

We have,

$$ \int_{B(y_0,l)} {u_{\epsilon_i}}^q=\int_0^l \int_{\partial B(y_0,r)} {u_{\epsilon_i}}^q(r\sigma_r)d\sigma_r dr \to 0, $$

We set, $ s_{i,q}(r)=\int_{\partial B(y_0,r)}{u_{\epsilon_i}}^q(r\sigma_r) $. Then,

$$ \int_0^l s_{i,q}(r)dr \to 0, $$

We can extract of, $ s_{i,q} $, a subsequence which noted $ s_{i,q} $ and which tends to $ 0 $ almost every-where on $ [0,l] $.

\smallskip

First, we choose, $ q_1=\dfrac{q(n+2)}{n-2} $ with $ q > \dfrac{n}{2} $, after we choose $ l_2 > 0 $, such that, $ \int_{B(y_0,l_2)} {u_{\epsilon_i}}^{q_1} \to 0 $. Finaly, we take $ l_1 \in ]0,l_2] $, such that, $  s_{i,q_1}(l_1)\to 0 $. We take $ l_0=\dfrac{l_1}{2} = l' $  in $ (*****) $ and $ l=l_1 $ in $ (*****) $, we obtain (if we use H\"older inequality for the two integrals of $ (*****) $),

$$ \exists \,\, l_0>0, \,\, \sup_{B(y_0,l_0)} u_{\epsilon_i} \to 0.$$

But, $ x_{\epsilon_i} \to y_0 $, for $ i $ large, $ x_{\epsilon_i} \in B(y_0,l_0) $, which imply,

$$ u_{\epsilon_i}(x_{\epsilon_i})=\max_{\Omega}u_{\epsilon_i} \to 0.$$

But if we write,

$$ u_{\epsilon_i}(x_{\epsilon_i})=\int_{\Omega} G(x_{\epsilon_i},y){u_{\epsilon_i}}^{N-1-\epsilon_i}(y)dy ,$$

we obtain,

$$ \max_{\Omega} u_{\epsilon_i}=u_{\epsilon_i}(x_{\epsilon_i})\leq (\sup_{\Omega} u_{\epsilon_i})^{N-1-\epsilon_i}\int_{\Omega} G(x_{\epsilon_i},y) dy=[u_{\epsilon_i}(x_{\epsilon_i})]^{N-1-\epsilon_i} w(x_{\epsilon_i}) ,$$

and finaly,

$$ 1 \leq u_{\epsilon_i}(x_{\epsilon_i})]^{N-2-\epsilon_i} w(x_{\epsilon_i}) . $$

But, $ w>0 $ on $ \Omega $, $ ||w||_{\infty} >0 $ and $ N-2-\epsilon_i >\dfrac{2}{n-2} $, we have,

$$ u_{\epsilon_i}(x_{\epsilon_i}) \geq \dfrac{1}{[{||w||_{\infty}}^{1/(N-2-\epsilon_i)}]} \geq  c_4(n,\Omega)>0. $$

It is a contradiction.

\smallskip

For the Theorem 3, we obtain a contradiction if we write:

$$ \max_{\Omega} u_{\epsilon_i} \leq (\max_{\Omega} u_{\epsilon_i})^{N-1}||w||_{\infty}+\max_{\Omega} u_{\epsilon_i} \epsilon_i \sup_{\Omega} \int_{\Omega} G(x,y)dy \leq ({\max_{\Omega} u_{\epsilon_i}})^{N-1}||w||_{\infty}+\alpha \max_{\Omega} u_{\epsilon_i}, $$

and finaly,

$$ \max_{\Omega} u_{\epsilon_i} \geq \left ( \dfrac{1-\alpha}{||w||_{\infty} }\right )^{1/(N-2)} .$$

\underbar { \bf Method 2:} proof of theorem 2 directly.

\bigskip

Suppose that:

$$ \sup_{\Omega} \times \inf_K u_i \to 0, $$

then, for $ \delta >0 $ small enough, we have:

$$ \sup_{\Omega} u_i \times \inf_{\{x,d(x,\partial \Omega )\geq \delta \}} u_i \to 0. $$

Like in the first method (see [H]), for $ \delta >0 $ small,

$$ \sup_{\{x,d(x,\partial \Omega) \geq \delta \}} u_i \leq M=M(n,\Omega). $$

We have,

$$ u_i(x)=\int_{\Omega} G(x,y)u_i^{N-1-\epsilon_i} dy. $$

Let us consider $ K' $ another compact of $ \Omega $, using maximum principle, we obtain:

$$ \exists \, c_1=c_1(K,K',n, \Omega)>0, \,\, {\rm such \,\, that } \,\, G(x,y) \geq c_1 \,\, \forall \,\, x\in K, \, y\in K, $$

thus,

$$ \inf_K u_i=u_i(x_i) \geq c_1 \int_{K'} u_i^{N-1-\epsilon_i} dy. $$

We take, $ K'=K_{\delta}=\{x,d(x,\partial \Omega)\geq \delta \} $, there exists $ c_2=c_2(\delta,n,K,\Omega) >0 $ such that:

$$ \sup_{\Omega} u_i \times \inf_K u_i \geq c_2 \int_{ K_{\delta} } u_i^{N-\epsilon_i} dy, $$

we deduce,

$$ ||u_i||_{N-\epsilon_i}^{N-\epsilon_i} \geq c_2' \sup_{\Omega} u_i \times \inf_{\{x,d(x,\partial \Omega)\geq \delta \}} u_i+ mes (\{x,d(x,\partial \Omega)\leq \delta \})M^{N-\epsilon_i}. $$

If we take $ \delta $ small and for $ i $ large, we have:

$$ ||u_i||_{N-\epsilon_i} \to 0 . $$

Now, we use the Sobolev imbedding, $ H_0^1 $ in $ L^N $, we multiply the equation of $ u_i $ by $ u_i $, we intgrate by part and finaly, by H\"older inequality, we obtain:

$$ \bar K_1 ||u_i||8{N-\epsilon_i}^2 \leq \bar K_2 ||u_i||_N^2 \leq \int_{\Omega} |\nabla u_i|^2=\int_{\Omega} u_i^{N-\epsilon_i}=||u_i||_{N-\epsilon_i}^{N-\epsilon_i}, $$

we know that, $ 0 < \epsilon_i \leq \dfrac{2}{n-2} $, the previous inequality:

$$ ||u_i||_{N-\epsilon_i} \geq \bar K_3 >0, \,\, \forall \,\, i,$$

it is a contradiction.

\newpage

{\small {\underbar {\bf References:}}

\bigskip

[AP] F.Atkinson, L.Peletier. Elliptic Equations with Nearly Critical Growth, J.Diff.Eq. vol 70, 1987, pp.349-365.

\smallskip

[A] T.Aubin. Some nonlinear problems in Riemannian Geometry. Springer-Verlag, 1998.

\smallskip

[B] S.S. Bahoura. Diff\'erentes Estimations du $ \sup u \times \inf u $ pour l'\'Equation de la Courbure Sclaire Prescrite en dimension $ n\geq 3 $. J. Math. Pures . Appl. (9)82 (2003), no.1, 43-66.

\smallskip

[B1] S.S. Bahoura. In\'egalit\'e de Harnack pour les solutions d'\'equations du type courbure scalaire prescrite. C.R.Math.Acad.Sci.Paris 341 (2005), no 1, 25-28.

\smallskip

[B M] H. Brezis, F. Merle. Uniform estimates and Blow-up Behavior for Solutions of $ -\Delta u=V(x) e^u $ in two dimension. Commun. in Partial Differential Equations, 16 (8 and 9), 1223-1253(1991).

\smallskip

[B N] H. Brézis, L. Nirenberg. Positive solutions of nonlinear elliptic equations invlving critical Sobolev exponents. Comm. Pure. Appl. Math. 36 (1983) (4) pp. 437-477.

\smallskip

[B P] H.Brezis, L.Peletier. Asymptotics for Elliptic Equations Involving Critical Growth. Partail Differential Equations and The Calculus of Variation. Vol 1 149-192, Progr. Nonlinear Differential Equations Appl., Birkhauser, Boston, Ma 1989.

\smallskip

[DJLW] W. Ding, J. Jost, J. Li, G. Wang. The Differential Equation $\Delta u=8\pi-8\pi h e^u $ On a Compact Riemann Surface. Asian. J. Math. 1, 230-248(1997).

\smallskip

[DLN] D.G. De Figueiredo, P.L. Lions, R.D. Nussbaum. A priori Estimates and Existence of Positive Solutions of Semilinear Elliptic Equations, J. Math. Pures et Appl., vol 61, 1982, pp.41-63.

\smallskip

[GNN] B. Gidas, W. Ni, L. Nirenberg. Symmetry and Related Propreties via the Maximum Principle, Comm. Math. Phys., vol 68, 1979, pp. 209-243.

\smallskip

[F] L. Fontana. Sharp bordeline Sobolev inequalities on compact Riemannian manifolds. Comment. Math. Helv. 68 (1993), no 3. 415-454.

\smallskip

[H] Z-C. Han. Assymptotic Approach to singular solutions for Nonlinear Elleptic Equations Involving Critical Sobolev Exponent. Ann. Inst. Henri Poincar\'e. Analyse Non-lin\'eaire. 8(1991) 159-174.

\smallskip

[L] YY.Li. Harnack type Inequality, the Methode of Moving Planes. Commun. Math. Phys. 200 421-444.(1999).

\smallskip

[O] E. Onofri. On the Positivity of the Effective Action in a Th\'eory of Random Surfaces. Commun. Math. Phys. 86, 321-326 (1982).

\smallskip

[S] J. Spruck. The Elliptic Sinh Gordon Equation and The Construction of Toroidal Soap Bubbles. Lectures Notes in Math., vol. 1340. Calculus of Variations and Partial Differential Equations (eds. Hildebrant, S., Kinderlehrer, D., and Miranda, C.), Springer, Berlin-Heidelberg-New York-Tokyo (1988).

\smallskip

[SN 1] T. Suzuki, K. Nagasaki. Asymptotic Behavior of Solution for an Elliptic Boundary Value Problem with Exponential Nonlinearity. Proc. Japan Acad. Ser. A Math. Sci. 65 (1989), no.3, 74-76.

\smallskip

[SN 2] T. Suzuki, K. Nagasaki. Asymptotic Analysis for two-Dimentionnal Elliptic Eigenvalue Problems with Exponentially Dominated Nonlinearities. Asymptotic Analysis 3 (1990) 173-188 }

\end{document}